\documentclass{amsart}                   
\usepackage{amsfonts}
\usepackage{amsmath}
\usepackage{amscd}
\usepackage{amssymb}
\usepackage[dvips]{graphicx}
\usepackage[dvips]{epsfig}

\def\Fin{\hfill$\blacksquare$\vskip2mm}

\def\mb#1{\mathbb{#1}}
\def\mc#1{\mathcal{#1}}

\def\R#1{\mb{R}^{^{#1}}}

\def\L{\mathcal{L}}
\def\X{\widetilde{X}}
\def\U{\mathcal{U}}

\def\normal{T_{_4}}

\def\tube#1{T_{_{#1}}}

\def\sub#1{_{_{#1}}}
\def\sup#1{^{^{#1}}}

\def\strat#1{\mathcal{S}\sub{#1}}           
\def\STRAT#1{\left({#1},\strat{#1}\right)}
\def\Min{\Sigma\sub{\text{min}}}
\def\lift#1{\widetilde{#1}}   
\def\X{\lift{X}}

\def\ARROW#1{
\text{\begin{picture}(35,18)(15,15)         
            \put(30,27){$_{#1}$}
            \put(18,20){\vector(1,0){30}}
\end{picture}}}

\def\LARROW#1{
\text{\begin{picture}(35,18)(15,15)         
            \put(30,27){$_{#1}$}
            \put(48,20){\vector(-1,0){30}}
\end{picture}}}

\def\Morf#1#2#3{{#1}{\ARROW{#2}}{#3}}

\newcounter{numero}
\newcommand{\Numero}{\setcounter{numero}{1}(\arabic{numero}) }
\newcommand{\numero}{\addtocounter{numero}{1}(\arabic{numero}) }

\newcounter{letra}
\newcommand{\Letra}{\medskip \setcounter{letra}{1}(\alph{letra}) }
\newcommand{\letra}{\medskip \addtocounter{letra}{1}(\alph{letra}) }

\newcounter{romnumero}

\newcounter{bibnumero}

\theoremstyle{plain}
\newtheorem{Teo}{Theorem}[subsection]
\newtheorem{Coro}[Teo]{Corollary}
\newtheorem{Lema}[Teo]{Lemma}
\newtheorem{Pro}[Teo]{Proposition}

\newtheorem{Con}[Teo]{Conjecture}

\theoremstyle{definition}
\newtheorem{Def}[Teo]{Definition}

\newtheorem{Ejms}[Teo]{Examples}
\newtheorem{Ejm}[Teo]{Example}

\theoremstyle{remark}                               
\newtheorem{Obs}[Teo]{Remark}

\newenvironment{Proof}{ [{\it Proof\/}]\rm\hskip2mm }{\hfill$\square$\vskip2mm}       
\newenvironment{sketch}{ [{\it Sketch of the Proof\/}]\rm\hskip2mm }{\hfill$\square$\vskip2mm}

\def\blema{\begin{Lema}}
\def\elema{\end{Lema}}

\def\bprop{\begin{Pro}}
\def\eprop{\end{Pro}}

\def\bteo{\begin{Teo}}
\def\eteo{\end{Teo}}

\def\bcor{\begin{Coro}}
\def\ecor{\end{Coro}}

\def\bconj{\begin{Con}}
\def\econj{\end{Con}}

\def\bdem{\begin{Proof}}
\def\edem{\end{Proof}}
\def\bsketch{\begin{sketch}}
\def\esketch{\end{sketch}}

\begin{document}

\title{Categorical properties of smooth unfoldings on stratified spaces}
\dedicatory{\small To Carely, Tomasito, Jolymar \&\ Santi.}

\author{T. Guardia}
\address{Universidad Central de Venezuela-Centro de Geometr\'ia-Escuela de Matem\'atica.}
\email{tomas.guardia@ciens.ucv.ve}

\author{G. Padilla}
\address{Departamento de Matem\'aticas, Universidad Nacional de Colombia, K30 con calle 45. 
Edificio 404, ofic. 315. Bogot\'a. (+571)3165000 ext 13166}
\email{gabrielpadillaleon@gmail.com}

\date{}
\keywords{Intersection Cohomology, Stratified Pseudomanifolds}
\subjclass{35S35; 55N33}

\maketitle

\begin{abstract} In a previous work we proved the uniqueness and functoriality 
of primary unfoldings on simple Thom-Mather 
spaces, which is a functor to the category of smooth manifolds.  
In this article we extend these results for any stratified Thom-Mather pseudomanifold with arbitary
finite length, through a new kind of intermediate desingularizations, the {\it unbendings}, which  
coincide with primary unfoldings in the simple case.

\end{abstract}

\section*{Introduction}

The intersection homology was defined by Goresky and MacPherson in order to extend 
the Poincar\'e duality to the family of spaces with singularities
\cite{gm1}. Among the earliest works concerning smooth desingularizations and their relation with
\linebreak intersection cohomology, we find \cite{ferrarotti, ferrarotti2, verona2,verona}. 
In \cite{brasselet} Brasselet, Hector and Saralegi defined the intersection cohomology
with differential forms on suitable smooth unfoldings and 
proved a stratified version of the De Rham theorem; 
the last author has continued a fruitful research in this direction \cite{illinois}.
Although the unfoldings are not uniquely determined, 
the intersection cohomology does not depend on their choice.   \vskip2mm

In \cite{dalmagro2} Dalmagro came back to
the geometrical point of view; he worked with primary unfoldings, a simpler
and slightly more restricted smooth desingularizations. In a previous article 
we proved the functorial behavior of primary unfoldings \cite{guardia},
so there is a canonical way to unfold simple Thom-Mather spaces.
In this article we extend these results for Thom-Mather stratified spaces with
arbitrary finite length, which is our first main result. This is accomplished as follows:
We construct a new kind of intermediate desingularizations, the {\it unbendings}. 
They are recursive steps which can be used in order to obtain smooth primary unfoldings,
and present nice functorial properties. The mutual incidence of tubular neighborhoods is avoided
since, in a Thom-Mather stratified space, any family of non-comparable strata
can be separated with a disjoint family of tubular neighborhoods.
The unbending of a simple Thom-Mather space coincides with its primary unfolding so, in a certain sense, unbendings are more general. Our second main result is that our unfolding is a functorial construction.
\vskip2mm 

This article has been organized as follows: Preliminary ideas are contained in \S1. In  
\S2 introduce the definition of unfoldings and unbendings. We devote \S3 to show the 
functorial properties of the unbending process. In \S4 we conclude with a proof of the
existence and functoriality of smooth unfoldings. Each time we use the word {\it manifold} we mean 
a smooth differentiable manifold of class $C^\infty$ without boundary.

\section{Stratified Pseudomanifolds}

\subsection{Stratified Spaces}\label{subsection stratified spaces}\label{def basic properties of stratified spaces}
\quad  In 1969 Thom \cite{thom} introduced
the notion of stratified spaces; they are metric spaces that can be decomposed in a locally finite disjoint union of smooth manifolds
satisfying a certain incidence condition. Let $X$ be a 2nd countable metric space. A family of subsets $\strat{}$ of $X$ is a {\bf stratification} 
if and only if $\strat{}$ is a locally finite partition of $X$ whose elements, with the induced topology,
are disjoint nonempty locally closed smooth
manifolds. A {\bf stratum} of $X$ is an element $S\in \strat{}$. Given
any other stratum $S'\in \strat{}$ we will say that $S'$ {\bf is incident
over} $S$ if $S\cap\overline{S'}\neq\emptyset$. The strata of $X$
are required to satisfy the following {\bf incidence condition}
\begin{equation}\label{eq estratificacion}
\text{If $S'$ is incident over $S$ then $S \subset\overline{S'}$}.
\end{equation}
If $\strat{}$ is a stratification of $X$ we say that $(X,\strat{})$ is a {\bf stratified space}, though 
we will and talk about {\it "a stratified space $X$"} whenever the choice of $\strat{}$ 
is clear in the context.\vskip2mm
For each stratified space $(X,\strat{})$ the following properties are straightforward \cite{pflaum}, 
\begin{enumerate}
	  \item The incidence condition is partial order relation.
	  \item There is at most a countable number of strata (i. e. $\strat{}$ is countable).	
	  \item For each stratum $S\in\strat{}$;\vskip2mm 
	  \begin{enumerate}
	  	\item $S$ is maximal (resp.  minimal) if and only if it is open (resp. closed).
  	  	\item The closure of $S$ is the union of the strata over which it is incident, 
		$\overline{S}=\underset{_{S'\leq S}}{\bigsqcup}\ S'$.
  	  	\item The set $U_S=\underset{_{S\leq S'}}{\bigsqcup}\ S'$ is open, we call it
          	the {\bf incidence neighborhood} of $S$.
	  \end{enumerate}
\end{enumerate}
\vskip2mm
A stratum $S\in\strat{}$ is {\bf regular} if it is open in $X$, otherwise we say that $S$ is
{\bf singular}. The {\bf singular part} (resp. {\bf regular part}) of $X$ is the union of
the singular (resp. regular) strata, which we note $\Sigma$ (resp. $X-\Sigma$). 
The \textbf{minimal part} is the union of closed (and therefore minimal) strata, denoted $\Min$.\vskip2mm 

A \textbf{stratified subspace} of $(X,\strat{})$ is a subset $Y\subset X$ such that
\[
	\strat{Y}=\{S\cap Y: S\in \strat{} \}
\] 
is a stratification of $Y$ with the induced topology.\Fin

\begin{Ejms}\label{ejem examples of stratified spaces}\hfill
\begin{enumerate}
  \item Each manifold $M$ is a stratified space with empty singular part $\Sigma=\emptyset$.
  \item If $M$ is a manifold and $(X,\strat{})$ is a stratified space
  then 
	\[
		\strat{M\times X}=\{M\times S : S\in\strat{} \}
	\]
	is a stratification of $M\times X$.
  \item Let $\left(L,\strat{L}\right)$ be a compact stratified space. The \textbf{open cone}
  of $L$ is the quotient space
  \[\label{def cono}
  	c(L)=\frac{L\times[0,1)}{\sim}
  \]
  where $(l,0)\sim (l',0)$ for every $l,l'\in L$. If $r=0$, we write $[l,r]$ for the equivalence class of a point
  $(l,r)$. The class of all points $(l,0)$ is the \textbf{vertex}
   of the cone and will be denoted as $v$. The stratification of $c(L)$
   is 
   \[
	\strat{c(L)}=\{v\}\sqcup\{S\times (0,1): S\in\strat{L} \}
   \]
   \item  A {\bf basic model} is a product of the form $M\times c(L)$ where $M$ is a manifold and $\STRAT{L}$ 
	is a compact stratified space.
\end{enumerate}
\end{Ejms}

Let us fix a stratified space $\left(X,\strat{}\right)$. 
Since $\strat{}$ is locally finite, each point $x\in
S$ has a neighborhood which intersects a finite number of strata.
We conclude that every strict incidence chain in $\strat{}$ is finite. This motivates
the next
\begin{Def}\label{def profundidad}
  \rm  The \textbf{length} of a stratum $S\in \strat{}$  is the largest integer $p\geq0$ such that there is a strict incidence chain 
\begin{equation}
   S=S\sub0< S\sub1<\cdots< S\sub{p}
\end{equation}
in $\strat{}$. The {\bf length} of $X$ is the supremum (possibly infinite) of the lengths of the strata.
We will denote it by $l(X)$. The \textbf{dimension} of $X$, denoted as $\dim(X)$, 
is defined in a similar way.\Fin
\end{Def}

\begin{Def}\label{def stratified morphism} 
A \textbf{stratified morphism} (resp. \textbf{isomorphism})
between two stratified spaces $\STRAT{X}$ and $\STRAT{Y}$ is a continuous function (resp. homeomorphism) 
$\Morf{X}{f}{Y}$
that sends smoothly (resp. diffeomorphically) each stratum of $X$  into a stratum of $Y$. 
A stratified morphism $f$ is an  \textbf{embedding} if $f(X)\subset Y$ is a stratified subspace and
$\Morf{X}{f}{f(X)}$ is an isomorphism.\Fin
\end{Def}
\begin{Obs}\label{obs morphism local writting}
Each morphism $\Morf{M\times c(L)}{f}{M'\times c(L')}$
 can be written as:
 \begin{equation}\label{def morfismo basico}
   f(u,[l,r])=(a_1(u,l,r),[a_2(u,l,r),a_3(u,l,r)])
 \end{equation}
 Where $a_1,a_2,a_3$ are maps defined on $M\times
 L\times[0,1)$ and are piecewise smooth, i.e, smooth on $M\times \{v\}$ and $M\times S\times (0,1)$ for each
 $S\in\strat{L}$.
\end{Obs}

\subsection{Stratified Pseudomanifolds}

A stratified pseudomanifold is a stratified \linebreak space together a family
of conic charts which reflect the way in which we approach
the singular part. The definition is given by induction 
on the length. 

\begin{Def} A 0-length {\bf stratified pseudomanifold} is a smooth manifold with the
trivial stratification. A stratified space $\left(X,\strat{}\right)$ with $l(X)>0$ 
is a {\bf stratified pseudomanifold} if, for each singular stratum $S$, 
\begin{enumerate}
	\item   There is a compact stratified pseudomanifold $\STRAT{L}$ with $l(L)<l(X)$. We call 
	$L$ the {\bf link} of $S$ because
	\item   Each point $x\in S$ has an open neighborhood $x\in U\subset S$
	and a stratified embedding $\Morf{U\times c(L)}{\alpha}{X}$ 
	on an open neighborhood of $x\in X$. 
\end{enumerate}
The image of $\alpha$ is called a {\bf basic neighborhood} of $x$. 
Notice that $\Im(\alpha)\cap S=U$. Without loss of generality, 
we assume that $\alpha(u,v)=u$ for each $u\in U$ (where $v$ is the vertex of $c(L)$, see \S\ref{ejem
examples of stratified spaces}-(4)). We summarize the above situation by saying that the pair
$(U,\alpha)$ is a {\bf chart} of $x$. The family of charts is an {\bf atlas} of $(X,\strat{})$. \Fin

\end{Def}

\begin{Ejms} \label{ejm examples of stratified pseudomanifolds}\hfill
\begin{enumerate}
  \item A basic model $U\times c(L)$  is a stratified pseudomanifold if $\STRAT{L}$ is a compact
  stratified pseudomanifold.
  \item If $M$ is a manifold and $X$  is a stratified pseudomanifold then
  $M\times X$ is a stratified pseudomanifold. 
  \item Every open subset of a stratified pseudomanifold is again a stratified \linebreak pseudomanifold. 
  \item Since algebraic manifolds satisfy the Withney's conditions, every algebraic manifold is a stratified
   pseudomanifold \cite{pflaum}. 
  \item The orbit space of a stratified pseudomanifold endowed with a suitable stratified action of a compact Lie
  group is again a stratified pseudomanifold \cite{gysin3,popper}.  
  \item The foliation space of a suitably controlled locally conic foliated manifold is a stratified pseudomanifold
  \cite{wolak, saralegi2}. 
  \item New examples of stratified pseudomanifolds are arising from the field of theoretical physics. See for instance
  \cite{huebschmann}.
\end{enumerate}
\end{Ejms}


\section{The process of removing singularities}
The main feature of any suitable desingularization 
is its cabapility of preserving (co)homological or geometrical information near the singular strata. In this article
we mention two kinds of desingularizations: smooth unfoldings 
\cite{brasselet,dalmagro2,davis,forme diff hector-sar,illinois} and the {\it unbendings} we are to introduce.
The main difference between the last \linebreak  two objects is that the unbending removes only the  
minimal part of the stratified pseudomanifold, while the unfolding removes completely the singular part.\vskip2mm
 
In \cite{guardia} we proved the equivalence between the Thom-Mather conditions and the existence of
unfoldings for simple spaces. In the following sections we will \linebreak  extend this result for any Thom-Mather
stratified pseudomanifold with arbitrary finite length.
In \cite{dalmagro2} Dalmagro works with transverse morphisms as an additional\linebreak  
requirement which we do not ask here. We will show that after a finite number of 
iterated compositions of unbendings we get Dalmagro's {\it "primary"} unfoldings.\vskip2mm

We fix a stratified pseudomanifold $\left(X,\strat{}\right)$ with finite length 
$l(X)=p<\infty$.

\subsection{Unfoldings}\label{def unfolding}
It is in terms of this geometric tool that 
the intersection \linebreak cohomology with smooth differential forms is defined. 
  An \textbf{unfolding} of $X$ consists of a manifold
  $\widetilde{X}$ and a continuous surjective proper map
  $\widetilde{X}\ARROW{\L}{X}$ satisfying the following conditions:\vskip2mm
  \begin{enumerate}
    \item \underline{Inductive condition:} 
	There is a family of unfoldings of the links of the singular strata
    $\{\L\sub{S}:\widetilde{L\sub{S}}\ARROW{}L\sub{S} : S\text{\ is singular }\}$.
    \item \underline{Regularity:} The restriction $\L\sup{-1}(X-\Sigma)\ARROW{\L} X-\Sigma$
    is a smooth finite trivial covering 
    (hence, a diffeomorphism on each copy of the regular part).
    \item \underline{Existence of unfoldable charts:} 
	For each point $z\in\L\sup{-1}(\Sigma)$, there is a \linebreak commutative diagram:
    \begin{center}
    \begin{picture}(70,70)(15,15)  \label{diagexplo}       
            \put(-50,10){$U\times c(L)$}
            \put(60,10){$X$}
            \put(-60,70){$U\times \widetilde{L}\times\mb{R}$}
            \put(60,70){$\X$}
 
            \put(-45,40){$^c$}
            \put(23,83){$_{\widetilde\alpha}$}
            \put(25,23){$_{\alpha}$}
            \put(70,40){$^{\L}$}
            \put(-35,65){{\vector(0,-1){40}}}
            \put(65,65){{\vector(0,-1){40}}}
            \put(5,75){\vector(1,0){40}}
            \put(5,15){\vector(1,0){40}}
        \end{picture}
    \end{center}\vskip2mm
    Such that:
    \begin{enumerate}
      \item $(U,\alpha)$ is a chart at $x=\L(z)$.
      \item $c(u,\widetilde{l},t)=(u,[\mathcal{L}_L(\widetilde{l}),|t|])$.
    \item $\widetilde{\alpha}$ is a diffeomorphism onto $\L^{-1}(\text{Im}(\alpha))$.
  \end{enumerate}
    The above diagram is an \textbf{unfoldable chart} on $x$.\vskip2mm
        A stratified pseudomanifold $X$ is \textbf{unfoldable} if it has an unfolding.
	\vskip2mm 
	\item An \textbf{unfoldable morphism} is a commutative diagram
    \begin{center}
    \begin{picture}(70,70)(15,15) \label{eq lifting of the unfolding}       
            \put(-5,10){$X$}
            \put(60,10){$X'$}
            \put(-5,70){$\X$}
            \put(60,70){$\widetilde{X'}$}
 
            \put(5,40){$^{\L}$}
            \put(30,83){$_{\widetilde{f}}$}
            \put(30,23){$_{f}$}
            \put(70,40){$^{\L'}$}
            \put(0,65){{\vector(0,-1){40}}}
            \put(65,65){{\vector(0,-1){40}}}
            \put(15,75){\vector(1,0){40}}
            \put(15,15){\vector(1,0){40}}
        \end{picture}
    \end{center}\vskip2mm
	such that the vertical arrows are unfoldings, $f$ is a stratified morphism and 
	$\widetilde{f}$ is a smooth map.\Fin
\end{enumerate}

\begin{Ejms}\label{ejems unfoldings}\hfill
\begin{enumerate}
    \item  Any diffeomorphism $M\ARROW{f}M$ is an unfolding, for any smooth manifold $M$ considered as a 0-length
	stratified pseudomanifold.
    \item The map $c$ given in \S\ref{def unfolding}-(2) is an unfolding of the
     basic model $U\times c(L)$.
    \item If $\Morf{\X}{\L}{X}$ is an unfolding then\vskip2mm 
	\begin{enumerate}
	\item For any open subset $A\subset X$ the restriction $\L^{^{-1}}(A)\ARROW{\L} A$ is an unfolding.
	\item For any singular stratum $S$ the restriction $\L^{^{-1}}(S)\ARROW{\L} S$ is a smooth bundle with
    	typical fiber $F=\widetilde{L}$ the unfolding of the respective link.
	\end{enumerate}\vskip2mm
    \item If $M\supset \Sigma$ is a singular manifold then $\Sigma$ has a smooth tubular neighborhood $T$.
    The family $\strat{M}=\{\Sigma,M-\Sigma\}$ is a stratification and $\STRAT{M}$ is a stratified pseudomanifold.
    If we take off $T$ and substitute it with the smooth fiber bundle $\widetilde{T}$ induced 
    by the radial homotetia, we find an \linebreak unfolding $\widetilde{M}\ARROW{\L} M$. For more details see 
    \cite{bredon,davis}. 
    \item An iteration of the above example shows that any smooth manifold $M$ endowed with a smooth Thom-Mather
    stratification has an unfolding \cite{brasselet}.
    \item Given a smooth effective action $G\times M\ARROW{\Phi}M$ of a compact Lie group $G$ on a smooth manifold $M$,
    the partition $\strat{M}$  of $M$ in orbit types endows $(M,\strat{M})$ with a smooth Thom-Mather stratification. 
    There is an unfolding $\widetilde{M}\ARROW{\L}M$, such that $\widetilde{M}$
    is endowed with the unique smooth free action of $G$ such that $\L$ is $G$-equivariant. 
    This structure passes in a canonical way
    to the respective orbit spaces $B=M/G$ and $\widetilde{B}=\widetilde{M}/G$ so that the orbit map
    $M\ARROW{\pi}B$ is unfoldable morphism \cite{forme diff hector-sar}. A completely analogous situation can be given
    for any stratified pseudomanifold that supports a suitable stratified action \cite{gysin3}.
\end{enumerate}
\end{Ejms}

\subsection{Unbendings}\label{def unbending}
The unbending of a stratified pseudomanifold is again a \linebreak stratified
pseudomanifold. In this section we study how the unfolding relates
with a finite composition of unbendings (as many as the depth of the
total space). \vskip2mm

An \textbf{unbending} of
  $X$ consists of a stratified pseudomanifold $\widehat{X}$ satisfying 
  $l\left(\widehat{X}\right)=l(X)-1$; and  continuous surjective proper map
  $\widehat{X}\ARROW{\widehat{\L}}X$ such that:
\begin{enumerate}
    \item The restriction $\widehat{\L}\sup{\hskip1mm -1}\left(X-\Min\right)
   \ARROW{\widehat{\L}}\left(X-\Min\right)$
    is a stratified double  \linebreak covering. The subset 
    $\widehat{\L}\sup{\hskip1mm -1}\left(X-\Min\right)$ is an open dense in $\widehat{X}$, 
    and it is the union of
    two disjoint isomorphic copies of $\left(X-\Min\right)$ which we
    denote $\left(X-\Min\right)\sup{\pm}$ and the restriction of $\widehat{\L}$
    to each of these copies is a stratified isomorphism.\vskip2mm
    \item For each $z\in\widehat{\L}^{-1}(\Min)$, there is an commutative diagram:
    \begin{center}
    \begin{picture}(70,70)(15,15)  \label{diagunbending}       
            \put(-50,10){$U\times c(L)$}
            \put(60,10){$X$}
            \put(-60,70){$U\times L\times\mb{R}$}
            \put(60,70){$\widehat{X}$}
 
            \put(-45,40){$^{\widehat{c}}$}
            \put(23,83){$_{\widehat\alpha}$}
            \put(25,23){$_{\alpha}$}
            \put(70,40){$^{\widehat{\L}}$}
            \put(-35,65){{\vector(0,-1){40}}}
            \put(65,65){{\vector(0,-1){40}}}
            \put(5,75){\vector(1,0){40}}
            \put(5,15){\vector(1,0){40}}
        \end{picture}
    \end{center}\vskip2mm
    Such that:
    \begin{enumerate}
      \item $(U,\alpha)$ is a chart of $x=\widehat{\L}(z)$.
      \item
      $\widehat{c}(u,l,t)=(u,[l,|t|])$.
    \item $\widehat{\alpha}$ is a stratified isomorphism on
    $\L^{^{-1}}(\text{Im}(\alpha))$.
  \end{enumerate}
      The above diagram is an \textbf{unbendable chart} on $x$.\vskip2mm 
      A stratified pseudomanifold $X$ is \textbf{unbendable}
     if it has an unbending. \vskip2mm
    \item An \textbf{unbendable morphism} is a commutative square diagram
    \begin{center}
    \begin{picture}(70,70)(15,15) \label{eq lifting of the unbending}       
            \put(-5,10){$X$}
            \put(60,10){$X'$}
            \put(-5,70){$\widehat{X}$}
            \put(60,70){$\widehat{X'}$}
 
            \put(5,40){$^{\widehat{\L}}$}
            \put(30,83){$_{\widehat{f}}$}
            \put(30,23){$_{f}$}
            \put(70,40){$^{\widehat{\L'}}$}
            \put(0,65){{\vector(0,-1){40}}}
            \put(65,65){{\vector(0,-1){40}}}
            \put(15,75){\vector(1,0){40}}
            \put(15,15){\vector(1,0){40}}
        \end{picture}
    \end{center}\vskip2mm
    such that the vertical arrows are unbendings and the horizontal arrows are stratified morphisms. \Fin
\end{enumerate}

\begin{Obs}\label{obs unfoldable implies unbendable}
  Every unfoldable stratified pseudomanifold is unbendable. 
It is enough to take $\widehat{X}$ as the closure of the quotient space of
$\widetilde{X}-\L^{-1}(\Sigma-\Sigma_{\text{min}})$ with the
following equivalence relation $z\sim z'$ if $\L(z)=\L(z')$. On the other hand,
if $l(X)=1$ then the links are compact smooth manifold, so the unbending is an unfolding. We
leave the details to the reader.
\end{Obs}

\begin{Ejms}\label{ejems unbendings}\hfill
\begin{enumerate}
    \item  Any diffeomorphism $M\ARROW{f}M$ is an unbending, for any smooth manifold $M$ considered as a 0-length
	stratified pseudomanifold.
    \item The map $\widehat{c}$ given in \S\ref{def unbending}-(2).(b) is an unbending of the
     basic model $U\times c(L)$.
    \item If $\Morf{\widehat{X}}{\widehat{\L}}{X}$ is an unbending then 
    \begin{enumerate}
	\item For any open subset $A\subset X$ the restriction $\widehat{\L}^{^{-1}}(A)\ARROW{\widehat{\L}} A$ 
	is an unbending.
	\item For any minimal stratum $S$ the restriction $\widehat{\L}^{^{-1}}(S)\ARROW{\widehat{\L}} S$ 
	is a \linebreak  stratified bundle with typical fiber $F=L$.
    \end{enumerate}    
\end{enumerate}
\end{Ejms}

\section{The unbending of a Thom-Mather pseudomanifold}
Tubular neighborhoods arised in riemannian geometry a useful tool for\linebreak approaching to 
closed submanifold in a controlled way; 
they were given by means of a riemannian metric and
a smooth transverse section of the closed submanifold called a {\it "slice"}; see
for instance \cite{bredon}. A generalization for stratified pseudomanifolds was first given by Thom in 
his historical article, while there are geometrical versions dealing with stratified slices 
\cite{popper,thom}.

\subsection{Thom-Mather pseudomanifolds}
A tubular neighborhood $T$ around a \linebreak singular  stratum  $S$ is a stratified
fiber bundle with two main features, a conic fiber and a global tubular radium. 
A Thom-Mather pseudomanifold is a stratified \linebreak pseudomanifold such that each singular stratum 
is contained in such a \linebreak neighborhood. \vskip2mm

We fix in the sequel a stratified pseudomanifold $(X,\strat{})$.

\begin{Def}\label{def tubes}\label{obs global radium} 
     Given a singular stratum $S$ in $X$,
     a \textbf{tubular neighborhood} on $S$  is a fiber bundle $\xi=(T,\tau,S,c(L))$ satisfying
    \begin{itemize}
        \item[\Numero] $T$ is an open neighborhood of $S$ in $X$.
        \item[\numero] The  fiber is $c(L)$, the open cone of the link of $S$.
        \item[\numero] The inclusion $S \subset T$ is a section: $\tau(x)=x$ for any $x\in S$.
        \item[\numero] The structure group is contained in Iso$\STRAT{L}$. If $(U,\alpha)$, $(V,\beta)$
        are two bundle charts and $U\cap V\neq\emptyset$ then the cocycle is
        {\small\[
            (U\cap V)\times c(L)\ARROW{\hskip-3mm\beta^{^{-1}}\alpha}
            (U\cap V)\times c(L)\hskip1cm   
	    \beta^{^{-1}}(\alpha(u,[l,r]))=(u,[g\sub{\alpha\beta}(u)(l),r])
        \]}
        where $g\sub{\alpha\beta}(u)$ is a stratified isomorphism of $\STRAT{L}$ for all $u\in U\cap V$.
    \end{itemize}
        Given a tubular neighborhood $T\ARROW{\tau} S$ by \S\ref{def tubes}-(4) 
	the cocycles of the fiber bundle preserve the conical radium, so it can be extended to a  
        well defined {\bf tubular radium} $T\ARROW{\rho} [0,\infty)$ which, locally, is given by 
	$\rho(\alpha(u,[l,r]))=r$ in the image of a bundle chart $(U,\alpha)$. Notice
        that $\rho^{^{-1}}(\{0\})=S$ and $\rho^{^{-1}}\left(\R{+}\right)=(T-S)$. There is also an action 
	$\R{+}\times T\ARROW{}T$, the {\bf radial stretching}, given by 
	$\lambda\cdot\alpha(u,[l,r])=\alpha(u,[l,\lambda r])$.\vskip2mm

        We will say that $X$  is {\bf Thom-Mather} provided that 
	each singular stratum $S$ has a tubular neighborhood $T_{_S}$.
	\Fin
\end{Def}

We state now an easy and useful result we will use hereafter.

\blema\label{lema non-comparable strata 2}
	In a Thom-Mather pseudomanifold, any family of non-comparable strata can be 
	separated with a family of disjoint tubular neighborhoods. 
\elema

\bdem
    Because...
    \begin{itemize}
	\item[$\bullet$] \underline{Tubular neighborhoods can be streched}: 
	Take a tubular neighborhood $T$ on a singular stratum $S$ and any other open neighborhood $O\supset S$.
	It is possible to find a smooth non-negative function so that the tubular radium $\rho$ of 
	$T$ is suitably streched in order to obtain a smaller tubular neighborhood $T'$ satisfying
	$S\subset T'\subset O$, with the same procedure employed by \cite{bredon} for smooth tubes.  
    \item[$\bullet$] \underline{Non-comparable strata can be separated in any stratified space:}
	Notice that\vskip2mm
	\begin{itemize}
		\item[\Letra] \underline{It is enough to show it for minimal strata:}
		If $\mc{F}\subset\strat{}$ is a family of non-comparable strata, take
		the union of the incidence neighborhoods
		$Z=\underset{_{S\in\mc{F}}}{\cup} U\sub{S}$. Then $Z$ is open in $X$,
		therefore $S\in\mc{F}$ iff $S$ is a minimal stratum in $Z$. 
		Since $Z$ is open, it is enough to give a family of disjoint neighborhoods
		in $Z$ separating the strata in $\mc{F}$.
		\item[\letra] \underline{ Minimal strata can be separated by disjoint open subsets}: 
		Any two different minimal strata in $X$ are disjoint closed subsets, that can be separated
		with two disjoint open subsets because $X$ is $\normal$.	
		The whole family of minimal strata can be separated because of \S\ref{subsection stratified spaces}(1), 
		(3)-(b), and the facts that
		$X$ is $\normal$, and $\strat{}$ is locally finite \cite{tesis}.
	\end{itemize}
    \end{itemize}
\edem

\begin{Obs}
	Lemma \S\ref{lema non-comparable strata 2} implies that we can now give a family of tubular neighborhoods
	for the singular strata, $\{\tau\sub{S}:T\sub{S}\ARROW{}S:S\text{\ is singular }\}$ such that non-comparable
	strata have disjoint tubes. If two tubular neighborhoods \linebreak $\tube{S},\tube{R}$ 
	have non-empty intersection; then the corresponding strata $S,R$ are \linebreak comparable:
	\[
		\tube{S}\cap \tube{R}\neq\emptyset\hskip2mm\Rightarrow\hskip2mm
		S\leq R\hskip2mm\text{or}\hskip2mm R\leq S 
	\]
	This incidence condition was quoted by Mather \cite[pp.43-44]{mather}.
\end{Obs}

\subsection{Local simplifications}\label{obs disjoint tubes}\label{subsection non-comparable strata}
We briefly describe the unbending process that we will make explicit in section \S\ref{subsection unbending}.\vskip2mm 

Take a family of disjoint tubes $\left\{T\sub{S}:S \text{\ is minimal\ } \right\}$ for the minimal strata. \linebreak
Unbend separately each tube obtaining a canonical map $\widehat{T\sub{S}}\ARROW{\hskip-2mm\widehat{\L}\sub{S}}T\sub{S}$
for each minimal stratum $S$. Each $\widehat{T\sub{S}}$ is a unique stratified fiber bundle 
$\widehat{T\sub{S}}\ARROW{\hskip-2mm\widehat{\tau\sub{S}}}S$ with fiber $L\sub{S}\times\R{}$, where $L\sub{S}$ is the link of $S$. There is an unbending of the tubular radium $\widehat{T\sub{S}}\ARROW{\hskip-2mm\widehat{\rho\sub{S}}}\R{}$ satisfying 
\[
	\left|\widehat{\rho\sub{S}}(x)\right|=\rho\sub{S}\left(\widehat{\L}\sub{S}(x)\right)
	\hskip2cm
	x\in \widehat{T\sub{S}}
\]
so $S\subset\widehat{T\sub{S}}$ and the inclusion is the 0-section. They also satisfy
$l\left(T\sub{S}\right)=p-1$. The difference $\widehat{T\sub{S}}\sup{*}=\left(\widehat{T\sub{S}}-S\right)$ has two connected components, say $\widehat{T\sub{S}}\sup{\pm}$, which are again stratified fiber bundles over $S$ with respective fibers $L\sub{S}\times\R{\pm}$. Define the global unbending of $X$ as the stratified amalgamation of two copies of $\left(X-\Min\right)$, say $\left(X-\Min\right)\sup{\pm}$, and the disjoint union of the unbended tubes. Again $l\left(\widehat{X}\right)=p-1$ so the unbending process decreases the length of the total space. \vskip2mm

Notice that in the process of unbending we do not touch the intermediate strata, but only disjoint tubes over non-comparable (minimal) strata. Hence, although the statements will remain as general as possible; in the context of the proofs and without loss of generality, by \S\ref{lema non-comparable strata 2}, we will make some or even all of the following assumptions:
\begin{enumerate}
	\item $X$ is a connected stratified pseudomanifold.
	\item $X$ has finite length $l(X)=p<\infty$.
	\item All strata in $X$ are comparable, i.e. $\strat{}$ is a well ordered set and there is a unique
	strict incidence chain $S\sub0<\cdots<S\sub{p}$. In particular, $S\sub0=\Min$.
	\item $X$ is Thom-Mather, and there is a family of tubular neighborhoods \linebreak
	$\{T\sub{k}\ARROW{\tau\sub{k}}S\sub{k}:0\leq k\leq p-1\}$
	such that $T\sub{k+1}\subset\left(T\sub{k}-S\sub{k}\right)$ for all $k$.
\end{enumerate}

\

\subsection{The unbending of a Thom-Mather pseudomanifold}\label{subsection unbending}
Now, we study the connection between the tubular neighborhoods and
the unbending. We follow now the main ideas of \cite{dalmagro2,davis}.

\begin{Lema}\label{lema unbending of a tube}
	Each tubular neighborhood $T\ARROW{\tau}S$ has an unbending
	$\widehat{T}\ARROW{\widehat{\L}}T$ such that: 
	\begin{enumerate}
		\item The composition $\widehat{T}\ARROW{\widehat{\tau}}S$ given by $\widehat{\tau}=\tau\widehat{\L}$
		is a stratified fiber bundle.
		\item The fiber of $\widehat{T}$ is $L\times\R{}$ where $L$ is the link of $S$.
		\item The cocycles of $\widehat{T}$ are the same of $T$.
		\item There is an unbending of tubular radium $\rho$, i.e. a continuous funtion 
		$\widehat{T}\ARROW{\widehat{\rho}}\R{}$ such that $\left|\widehat{\rho}(x)\right|=\rho\left(\widehat{\L}
		(x)\right)$ for al $x\in\widehat{T}$.
	\end{enumerate}
\end{Lema}

\bdem
Let us fix a bundle atlas $\U=\{(U\sub\alpha,\alpha)\}\sub{\alpha\in\mathfrak{I}}$ for $T$.
\begin{itemize}
	 \item[\Letra] \underline{Unbending of a chart:} For any chart
	 $(U,\alpha)\in\U$, the unbending of $\tau^{-1}(U)$ is just the
 	composition 
 	\[
		U\times L\times\R{}\ARROW{\widehat{c}}U\times c(L)\ARROW{\alpha} \tau^{-1}(U)
 	\]
	where $\widehat{c}$ is the map defined in \S\ref{def unbending}-(2)-(b).
 	\item[\letra] \underline{Definition of the bundle $\widehat{T}$:}
 	Take the quotient space 
	{\small
 	\[
		\hskip1cm\widehat{T}=\frac{\bigsqcup_\alpha U_\alpha\times L\times
		\R{}}{\sim}\hskip1cm
		(u,l,t)\sim(u,g\sub{\alpha\beta}(u)(l),t)\ \forall
		\alpha,\beta\  \forall\ u\in U\sub\alpha \cap U\sub\beta
	\]}
	Denote by $[u,l,t]$ the equivalence class of $(u,l,t)$. 
	Following \cite[p.14]{steenrod} we get a unique fiber bundle
	\[
		\widehat{T}\ARROW{\widehat{\tau}}S\hskip2cm
		\widehat{\tau}([u,l,t])=u.
	\]
	with fiber $F=L\times\R{}$ and the same structure group of $T$.
	\item[\letra] \underline{Unbending of $T$:} Define 
	{\small\[
		\hskip1cm
		\widehat{T}\ARROW{\widehat{\L}}T\hskip5mm
		\widehat{\L}([u,l,t])=\alpha(u,[l,|t|])\hskip5mm 
		\forall l\in L,\ \forall t\in\R{}, \forall\ u\in U\sub\alpha,\
		\forall \alpha
	\]}
	Since the  cocycles $g\sub{\alpha\beta}$ are stratified isomorphisms of 
	$L$, we conclude that $\widehat{T}$ is a stratified
	pseudomanifold and $l\left(\widehat{T}\right)= l\left(T\right)-1=l(X)-1$. 
	In order to show that the above arrow is an unbending, let
	$(U,\alpha)\in\U$. Define
	\[
		U\sub\alpha\times L\times\R{}\ARROW{\widehat{\alpha}}\widehat{T}
		\hskip2cm
		\widehat{\alpha}(u,l,r)=[u,l,r].
	\]
	Then $\widehat{\alpha}$ is stratified, because it is the restriction of the quotient map. 
	The induced diagram \S\ref{def unbending}-(2) commutes $\forall\alpha$.\vskip2mm
	\item[\letra] \underline{Unbending of the tubular radium:} 
	The function
	{\small
	\[
		\widehat{T}\ARROW{\widehat{\rho}}\R{}\hskip5mm
		\widehat{\rho}\left(\widehat{\alpha}(u,l,t)\right)=t
		\hskip5mm
		\forall l\in L,\ \forall t\in\R{}, \forall\ u\in U\sub\alpha,\
		\forall \alpha
	\]}
 	trivially satisfies the required property.
\end{itemize}
\edem

\begin{Pro}\label{pro thomather implies unbendable}
	Let $X$ be a stratified pseudomanifold with finite length. 
	If every minimal stratum has a tubular neighborhood,
	then $X$ is unbendable.
\end{Pro}

\bdem
By \S\ref{subsection non-comparable strata} assume that $X$ has a unique minimal stratum $S\sub0=\Min$,
with a tubular neighborhood $T\sub0$. By \S\ref{lema unbending of a tube} and since the cocycles 
are radium-preserving, 
\[
		\widehat{\L}^{^{-1}}(T-S\sub0)=T\sub0\sup{+}\sqcup T\sub0\sup{-}
\]
has two connected components. They are stratified bundles
over $S\sub0$ with respective fibers $F\sup{\pm}=L\sub0\times\R{\pm}$. 
We obtain the unbending of the whole stratified pseudomanifold $X$  by taking two copies 		
$(X-S\sub0)\sup{\pm}$ of $X-S\sub0$; and suitably gluing them together along $\widehat{T}$.
In other words, we take 
\[
	\widehat{X}=\frac{(X-S\sub0)\sup{+}\sqcup\widehat{T\sub0}\sqcup(X-S\sub0)\sup{-}}{\sim}
\]
as the amalgamated sum by the inclusions $T\sub0\sup{\pm}\subset(X-S\sub0)\sup{\pm}$.
\edem

\begin{Coro}\label{cor unbendable does not change supositions}\label{pro thomather implies unbendable 2}
	Every Thom-Mather pseudomanifold is unbendable.
	The unbending of a connected Thom-Mather pseudomanifold of length $p$ is
	a connected Thom-Mather pseudomanifold of length $(p-1)$.
\end{Coro}
\bdem
	This is a consequence of \S\ref{subsection non-comparable strata}-(4), \S\ref{lema unbending of a tube}-(4)
	and \S\ref{pro thomather implies unbendable}; the unbending process
	does not affect the tubular neighborhoods of non-minimal strata.
\edem

\subsection{Thom-Mather morphisms}\label{def morfismo thomather}
	A {\bf tube-morphism} $\tube{S}\ARROW{f}\tube{R}$ between \linebreak
	stratified tubular neighborhoods 
	is a stratified morphism $f$ such that  
	\begin{enumerate}
		\item It commutes with the tubular radia, $\rho\sub{R}f=\rho\sub{S}$.
		\item It is a bundle-morphism, $\tau\sub{R}f=f\tau\sub{S}$.
	\end{enumerate}\vskip2mm
	Condition (1) implies that $\varphi(S)\subset R$ so (2) makes sense. Notice that 
	a tube-morphism $f$ commutes with the respective bundle cocycles of $\tube{S},\tube{R}$, 
	see \cite{steenrod}. \vskip2mm

	A {\bf Thom-Mather morphism} $X\ARROW{f}Y$ between Thom-Mather stratified pseudomanifolds,
	is a tube-preserving stratified morphism. In other words, $f$
	is a stratified morphism which sends tubes on tubes, 
	$f\left(\tube{S}\right)\subset\tube{f(S)}$ for each singular
	stratum $S$ of $X$; and $f$ is a tube-morphism on each tubular neighborhood.\Fin

\begin{Ejm}\label{Ejem tube isomorphism}
	If $X$ is a Thom-Mather stratified pseudomanifold and $X\ARROW{f}X$ is a stratified 
	isomorphism, then $f$ is a Thom-Mather isomorphism. For each
	\linebreak tubular neighborhood $T\ARROW{\tau}S$ it is enough to define 
	$T'\sub{S}=f(T\sub{S})$, $\tau'=f\tau$, and 
	$\rho'=f\rho$. Then $T'\ARROW{\tau'}S$ is a tubular neighborhood and
	$T\ARROW{f}T'$ is a tube-isomorphism. 	
\end{Ejm}

\subsection{Functoriality of unbendings}
In \cite{guardia} we proved that the primary unfoldings are representative in the
category of the unfoldable pseudomanifolds. In the more general context of this work
unbendings have the same representation property and, for simple pseudomanifolds they 
coincide with the primary unfoldings.

\begin{Lema}\label{lema continuous extension}
Let $X,X'$ be unbendable stratified pseudomanifolds with finite length,
$\Morf{\widehat{X}}{\widehat{\L}}{X}$ and
$\Morf{\widehat{X'}}{\widehat{\L'}}{X'}$ two  unbendings. 
Then, for each stratified morphism $\Morf{X}{f}{X'}$ there is a unique continous function
$\Morf{\widehat{X}}{\widehat{f}}{\widehat{X'}}$ such that
\begin{enumerate}
	\item $f\widehat{\L}=\widehat{\L}'\widehat{f}$ i.e. the diagram \S\ref{def unbending}-(3) commutes; 
	\item The restriction of $\widehat{f}$ to $\widehat{\L}\sup{\hskip1mm -1}\left(X-\Min\right)$ 
	is a stratified morphism. 
\end{enumerate}
\end{Lema}
\bdem  According to \S\ref{def unbending}-(1) the open dense
\[
	(X-\Min)\sup{\pm}=\widehat{\L}\sup{\hskip1mm -1}(X-\Min)\subset\widehat{X}
\] is the union of two isomorphic
copies of $(X-\Min)$ and the restriction of $\widehat{\L}$ to each of these copies is a stratified isomorphism.
A similar situation happens for $\widehat{X'}$ and we write $(X'-\Min')\sup{\pm}$ for the respective
copies of $(X'-\Min')$.
\begin{itemize}
	\item[\Letra]\underline{Definition of $\widehat{f}$ on $(X-\Min)\sup{\pm}$:}
	Then the inverse maps 
	\[
		\hskip1cm(X'-\Min')\ARROW{\hskip-2mm\widehat{\L'}\sup{\hskip1mm -1}} (X'-\Min')\sup{+}
		\hskip1cm
		(X'-\Min')\ARROW{\hskip-2mm\widehat{\L'}\sup{\hskip1mm -1}}(X'-\Min')\sup{-}
	\]
	are stratified isomorphisms. There are two ways in order to
	define the composition $\widehat{f}=\widehat{\L'}\sup{\hskip1mm -1}f\widehat{\L}$,
	depending on the copies of $(X'-\Min')$ where we take the inverse. Namely 
	we can take $\widehat{f}=\widehat{f}\sup{\hskip1mm +}$ as the {\it "sign-preserving"}
	lifting, i.e. 	
	{\small\[
		\widehat{f}\sup{\hskip1mm +}(X-\Min)\sup{+}\subset(X'-\Min')\sup{+}
		\hskip1cm
		\widehat{f}\sup{\hskip1mm +}(X-\Min)\sup{-}\subset(X'-\Min')\sup{-}
	\]}  
	and $\widehat{f}=\widehat{f}\sup{\hskip1mm -}$ as the other possible definition. In any case of these cases, 
	$\widehat{f}=\widehat{f}\sup{\hskip1mm \pm}$ 
	is a well-defined stratified morphism and satisfies \S\ref{lema continuous extension}.

	\item[\letra]\underline{ Continuous extension of $\widehat{f}$:} By the previous step,
	we have already defined a stratified morphism $\widehat{f}$  on $(X-\Min)\sup{\pm}$ satisfying
	\S\ref{lema continuous extension}. We define a global continuous extension of 
	$\widehat{f}$  as follows. \vskip4mm
	
	Assume by \S\ref{subsection non-comparable strata} that $\Min=S\sub0$ is a single minimal stratum.
	Let \linebreak $\widehat{z}\in\widehat{\L}\sup{-1}\left(S\sub0\right)$. We must define 
	$\widehat{f}(\widehat{z})$. For this sake let
	\[
		\{\widehat{z}\sub{n}\}\sub{n}\subset\widehat{\L}\sup{-1}\left(X-S\sub0\right)
	\]
	be any sequence converging to $\widehat{z}$. Since $\L,f$
	are continuous and $\widehat{\L},\widehat{\L}'$ is a continuous proper
	maps; by an argument of compactness and up to minor adjustments, we
	may suppose that the sequence $\{\widehat{f}(\widehat{z}\sub{n})\}\sub{n}$
	converges in $\widehat{X'}$. We define
	\[
		\widehat{f}(\widehat{z})=\lim\limits_{n\rightarrow\infty}
		\widehat{f}(\widehat{z}_n)
	\]
	If our limit-definition makes sense then it is also continuous; so
	next we will show that $\widehat{f}$ is well defined. Since the
	above definition is local, we first study the\vskip2mm

	\item[\letra]\underline{Lifting in terms of conics charts:}
	Assume that $X=M\times c(L)$ and
	$X'=M'\times c(L')$ are trivial basic models and their respective
	unbendings are the canonical ones - see \S\ref{ejems
	unbendings}. Then $f$ can be written as in \S\ref{obs morphism
	local writting}. The point
	\[
		\widehat{z}=(u,l,0)\in M\times L\times\{0\}\subset\widehat{X}=M\times L\times\R{}
	\]
	is the limit of a sequence
	\[
		\left\{\widehat{z}\sub{n}=\left(u\sub{n},l\sub{n},t\sub{n}\right)\right\}\sub{n}\subset M\times
		L\times(\R{}-\{0\})
	\]
	So the sequences $\left\{u\sub{n}\right\}\sub{n}$, $\left\{l\sub{n}\right\}\sub{n}$ and $\left\{t\sub{n}
	\right\}\sub{n}$
 	converge to $u$, $l$ and 0 respectively. Since $\widehat{\L}=\widehat{c}$, $\widehat{\L'}=\widehat{c'}$ and 
	$f$ are continuous maps, the sequence
	\[
		\hskip1cm w\sub{n}=f\left(\widehat{c}\left(\widehat{z}\sub{n}\right)\right)=
		\left(a\sub1\left(u\sub{n},l\sub{n},\left|t\sub{n}\right|\right),
		\left[a\sub2\left(u\sub{n},l\sub{n},\left|t\sub{n}\right|\right),
		a\sub3\left(u\sub{n},l\sub{n},\left|t\sub{n}\right|\right)\right]\right)
	\]
	converges to $w=f\left(\widehat{c}\left(\widehat{z}\right)\right)=(a\sub1(u,l,0), v)$. By
	the continuity of the functions $a\sub{j}$ for $j=1,2,3$ and up to minor
	adjustments on $a\sub{2}$ concerning the compactness arguments; we get
	that
	\[
		\widehat{w}\sub{n}=(a\sub1(u\sub{n},l\sub{n},|t\sub{n}|), a\sub2(u\sub{n},l\sub{n},|t\sub{n}|),\pm
		a\sub3(u\sub{n},l\sub{n},|t\sub{n}|))
	\]
	converges to 
	\[
		\widehat{w}=(a\sub1(u,l,0), a\sub2(u,l,0),0)\in M'\times L'\times\{0\}
	\]
	\item[\letra] \underline{The lifting is well defined:}
	From the continuity of the functions $a\sub{i}$, $i=1,2,3$; it follows that the
	element $\widehat{w}$ does not depend on the choice of a particular
	sequence $\{\widehat{z}\sub{n}\}\sub{n}$.
\end{itemize}
\edem

\begin{Def}\label{def almost unbending}
	The $\widehat{f}$ obtained at \S\ref{lema continuous extension}
	an \textbf{almost-unbending} of $f$.
\end{Def}

\begin{Pro}\label{pro local unbending properties}
	The almost-unbending 
	{\small
	\[
		U\times L\times\R{}\ARROW{\widehat{f}}
		U'\times L'\times\R{}
		\hskip1cm
		\widehat{f}(u,l,t)=(\widehat{a\sub1}(u,l,t),\widehat{a\sub2}(u,l,t),\widehat{a\sub3}(u,l,t))
	\]}
	of a stratified morphism between basic models
	{\small\[
		U\times c(L)\ARROW{f} U'\times c(L')\hskip1cm
		f(u,[l,r])=(a\sub1(u,l,r),[a\sub2(u,l,r),a\sub3(u,l,r)])
	\]}
	is an unbending as in \S\ref{def unbending}-(3), if and only if, for each stratum $S\in\strat{L}$,
	\begin{enumerate}
		\item[\Letra] $\widehat{a\sub3}(u,l,0)=a\sub3(u,l,0)=0$ for all $u\in U$, $l\in S$.
		\item[\letra] With respect to the coordinate $t\in\R{}$, 
		the functions $\widehat{a\sub{1}},\widehat{a\sub{2}}$ are even 			
		and $\widehat{a\sub3}$ is either odd or even. 
		\item[\Letra] $\widehat{a\sub{j}}$ is a smooth extension of  $a\sub{j}$ for all $j$.
	\end{enumerate}
\end{Pro}

\bdem 
	If $\widehat{f}$ is a unbending of $f$, then
	$fc=c'\widehat{f}$ where $c$ and $c'$ are canonical unbendings as in
	\S\ref{def unfolding}-(2). Checking both sides of this equality we
	get
	{\small\[
		f(c(u,l,t))=f(u,[l,|t|])=(a_1(u,l,|t|), [a_2(u,l,|t|),
		a_3(u,l,|t|)])
	\]}
	and
	{\small\[
		c'(\widehat{f}(u,l,t))=c'(\widehat{a}_1(u,l,t),
		\widehat{a}_2(u,l,t),\widehat{a}_3(u,l,t))= (\widehat{a}_1(u,l,t),
		[\widehat{a}_2(u,l,t), |\widehat{a}_3(u,l,t)|])
	\]}
	we conclude that
	{\small\[
		(a_1(u,l,|t|), [a_2(u,l,|t|), a_3(u,l,|t|)])=(\widehat{a}_1(u,l,t),
		[\widehat{a}_2(u,l,t), |\widehat{a}_3(u,l,t)|])
	\]}
	There are two cases; $t=0$ and $t\neq 0$, from which we get
	\S\ref{pro local unbending properties}.
\edem

\begin{Lema}\label{pro unbendable cocycles}
	The cocycles of any tubular neighborhood are unbendable.
\end{Lema}

\bdem 
	For  $f=\varphi=\beta^{-1}\alpha$ as in
	\S\ref{def tubes}-(4); the functions $\widehat{a\sub1}(u,l,t)=u$,
	$\widehat{a\sub2}(u,l,t)=g\sub{\alpha\beta}(u)(l)$ and $\widehat{a\sub{3}}(u,l,t)=t$ 
	satisfy the hypothesis of \S\ref{pro local unbending properties}.
\edem

\begin{Lema}\label{lema conmutativity of the unbending with cocycles}
Consider a diagram of stratified morphisms
    \begin{center}
    \begin{picture}(70,70)(15,15)  \label{diagunbending}       
            \put(-50,10){$U\times c(L)$}
            \put(60,10){$U'\times c(L')$}
            \put(-60,70){$U\times c(L)$}
            \put(60,70){$U'\times c(L')$}
             \put(-45,40){$^{\varphi}$}
            \put(23,83){$_{f}$}
            \put(25,23){$_{f'}$}
            \put(90,40){$^{\varphi'}$}
            \put(-35,65){{\vector(0,-1){40}}}
            \put(85,65){{\vector(0,-1){40}}}
            \put(5,75){\vector(1,0){40}}
            \put(5,15){\vector(1,0){40}}
        \end{picture}
    \end{center}\vskip2mm
such that $\varphi,\varphi'$ are as in \S\ref{def tubes}-(4). Then $f'\varphi=\varphi'f$
if and only if:
\begin{equation}\label{eq cocycle conditions of the lifting}
\begin{split}
  a\sub1(u,l,r)=a\sub1'(u,g(u)(l),r)\\
  g'(a\sub1(u,l,r))a\sub2(u,l,r)=a\sub2'(u,g(u)(l),r)\\
  a\sub3(u,l,r)=a\sub3'(u,g(u)(l),r)
\end{split}
\end{equation}
\end{Lema}
\bdem 
Write $\varphi(u,[l,r])=(u,g(u)(l),r)$ and $\varphi'(u',[l',r])=(u',g'(u')(l),r)$ where
$g(u),g'(u')$ are, respectively, isomorphisms on $L,L'$. This is a straightforward calculation.
\edem

\begin{Lema}\label{tubemorphisms are unbendable}
	Each tube-morphism is unbendable.
\end{Lema}
\bdem
	Let $\tube{S}\ARROW{f}\tube{R}$ be a tube-morphism. By \S\ref{lema continuous extension} we must show that
	the almost-unbending $\widehat{f}$ is stratified, i.e. that\vskip2mm
	$\bullet$ \underline{$\widehat{f}(\widehat{\L}\sup{\hskip1mm -1}(S))\subset 
	\widehat{\L}\sup{\hskip1mm -1}(R)$:}
	Since $\rho\sub{R}f=\rho\sub{S}$, we deduce that the unbended radia satisfy a similar property,
	$\widehat\rho\sub{R}\widehat{f}=\widehat\rho\sub{S}$.\vskip2mm
	$\bullet$ \underline{$\widehat{f}$ is stratified on $\widehat{\L}\sup{\hskip1mm -1}(S)$:}
	Since $f$ commutes with the respective cocycles of $\tube{S},\tube{R}$, 
	this is in fact a local matter. For any pair of bundle charts
	$(U,\alpha)$ of $\tube{S}$, and $(U',\alpha')$ of $\tube{R}$; the composition
	$h=\beta\sup{-1}f\alpha$ satisfies the requirements of \S\ref{pro local unbending properties}.
\edem

\begin{Teo}\label{def global thomather implies local thomather}\hfill
	\begin{enumerate}
		\item Each Thom-Mather morphism is unbendable.
		\item The unbending process is a functor in the category of Thom-Mather stratified
		pseudomanifolds.
		\item The unbending of a Thom-Mather stratified pseudomanifold
		is unique up to Thom-Mather isomorphisms.
	\end{enumerate}
\end{Teo}
\bdem 
	This is a  consequence of \S\ref{def morfismo thomather} and 
	\S\ref{tubemorphisms are unbendable}.
\edem

\section{Categorical properties of smooth desingularizations}
The aim of this section is to establish the sufficient and necessary conditions
for the existence of a smooth unfolding. We will also prove that the unfoldings have
similar functorial properties which they inherit from unbendings. \vskip2mm

\subsection{The primary unfolding of a stratified Thom-Mather pseudomanifold}
We start with an useful and easy result,
\begin{Lema}\label{lema unbendings are TM}
	The unbending map of a Thom-Mather stratified pseudomanifold
	is a Thom-Mather morphism.
\end{Lema}
\bdem
	In order to prove this we can assume the simpler geometrical conditions of \S\ref{subsection non-comparable strata}.
	By \S\ref{pro thomather implies unbendable 2}-(3) and (4), the tubular radium $\rho\sub0$ 
	of the tubular neighborhood $T\sub0$ on the minimal stratum $S\sub0$ is unbendable.
	Since the other tubular neighborhoods are open subsets of $(T\sub0-S\sub0)$
	and the preimage $\widehat{\L}\sup{\hskip1mm -1}(X-S\sub0)=(X-S\sub0)\sup{\pm}$ is the union
	of two isomorphic copies of $(X-S\sub0)$;  
	acording to \S\ref{def morfismo thomather} and \S\ref{Ejem tube isomorphism}, 
	$\widehat{\L}$ is tube-preserving on $\widehat{\L}\sup{\hskip1mm -1}(X-S\sub0)$. 
	We conclude that $\widehat{\L}$ is a Thom-Mather morphism. 
\edem

\begin{Teo}\label{Teo TM implies unfoldable}
Every finite-length Thom-Mather stratified pseudomanifold is unfoldable.
\end{Teo}
\bdem 
Let $p=l(X)$. We proceed in two steps.
\begin{itemize}
	\item[$\bullet$] \underline{Definition of the unfolding:}
	Corollary \S\ref{cor unbendable does not change supositions} implies that the unbending process does not
	change the assumptions \S\ref{subsection non-comparable strata}-(1) and (2). Assumption \S\ref{subsection 
	non-comparable strata}-(3) also holds because any intermediate (non-minimal) stratum in $X$ is locally detected 
	near $S\sub0$ at the
	link $L\sub0$. Since the unbending process does not touch the links we deduce that $\widehat{X}$ still has only
	one incidence chain. This allows us to continue an iterative unbending operation. \vskip2mm

	Denote the first unbending of $X$ by $\widehat{X}=X\sup1$ and the corresponding
	map by $\widehat{\L}=\widehat{\L}\sup1$. We obtain a finite sequence of unbendings
	{\small
	\[
		X\LARROW{\widehat{\L}\sup1}X\sup1 \LARROW{\widehat{\L}\sup2}\hskip3mm
		\cdots \hskip3mm\LARROW{\widehat{\L}\sup{p-1}}
		X\sup{p-1}\LARROW{\widehat{\L}\sup{p}}X\sup{p}
	\]}
	Notice that $l\left(X\sup{p}\right)=0$, so $X\sup{p}$ is a manifold. Take
	$\X=X\sup{p}$ and $\L=\widehat{\L}\sup{p}\cdots\widehat{\L}\sup1$. 
	\item[$\bullet$] \underline{Existence of unfoldable charts:} In order to show that
	$\X\ARROW{\L}X$ is an unfolding we check condition \S\ref{def unfolding}-(2). This can be done 
	by induction on $p$. For $p=0$ is trivial and for $p=1$ the minimal stratum $S\sub0=\Sigma$
	coincides with the singular part. The link $L\sub0$ is a compact smooth manifold and the unbending
	is the unfolding; this case has been treated in \cite{guardia}. We assume the inductive hypothesis
	so for any $k\leq p$ the statement holds. This implies that
	\begin{enumerate}
		\item For any singular stratum $S$ the respective link $L$ is unfoldable 
		in the described way.
		\item Since $l(X\sup1)=p-1$, the statement holds for $X\sup1$. 
	\end{enumerate}
	Take $\L'=\widehat{\L}\sup{p}\cdots\widehat{\L}\sup2$. By inductive hypothesis 	
	\[
		X\sup1\LARROW{\L'}\X
	\]
	is an unfolding of $X\sup1$. Consider the composition
	\[
		X\LARROW{\widehat{\L}\sup1}X\sup1 \LARROW{\L'}\X\hskip1cm
		\L=\widehat{\L}\sup1\L'
	\]
	Take a point $z\in\L\sup{-1}\left(S\sub0\right)$. Let us verify condition \S\ref{def unfolding}-(2),
	i. e. the existence of an unfoldable chart at $z$.\vskip3mm
 
	Take an unbendable chart as in \S\ref{def unbending}-(2) at $z'=\L'(z)\in X\sup1$
    \begin{center}
    \begin{picture}(70,70)(15,15)  \label{diagunbending}       
            \put(-50,10){$U\times c\left(L\sub0\right)$}
            \put(60,10){$X$}
            \put(-60,70){$U\times L\sub0\times\mb{R}$}
            \put(60,70){$\widehat{X}=X\sup1$}
 
            \put(-45,40){$^{\widehat{c}}$}
            \put(23,83){$_{\widehat\alpha}$}
            \put(25,23){$_{\alpha}$}
            \put(70,40){$^{\L\sup1}$}
            \put(-35,65){{\vector(0,-1){40}}}
            \put(65,65){{\vector(0,-1){40}}}
            \put(5,75){\vector(1,0){40}}
            \put(5,15){\vector(1,0){40}}
        \end{picture}
    \end{center}\vskip2mm
	Since $l\left(L\sub0\right)=p-1$, as we already remarked, by inductive argument 
	the link $L\sub0$ can be unfolded with a finite	sequence of $p-1$ unbendings,
	{\small
	\[
		L\sub0\LARROW{\eta\sup1}L\sub0\sup1 \LARROW{\eta\sup2}\hskip3mm
		\cdots \hskip3mm\LARROW{\eta\sup{p-2}}
		L\sub0\sup{p-2}\LARROW{\eta\sup{p-1}}L\sub0\sup{p-1}
	\]}
	By \S\ref{def global thomather implies local thomather} and \S\ref{lema unbendings are TM};
	these unbendings behave in a functorial way. Since $U\times L\sub0\times\R{}$
	is a stratified Thom-Mather pseudomanifold and $l\left(U\times L\sub0\times\R{}\right)=l\left(L\sub0\right)=p-1$;
	we deduce that the composition of the maps $\nu\sub{j}=id\sub{U}\times\eta\sup{j}\times id\sub{\R{}}$
	provides a finite sequence of $p-1$ unbendings
	{\small
	\[
		U\times L\sub0\times\R{}\LARROW{\nu\sub1}U\times L\sub0\sup1\times\R{} \LARROW{\nu\sub2}\hskip3mm
		\cdots \hskip3mm\LARROW{\nu\sub{p-2}}
		U\times L\sub0\sup{p-2}\times\R{}\LARROW{\nu\sub{p-1}}U\times L\sub0\sup{p-1}\times\R{}
	\]}
	Again by inductive hypothesis, we deduce that the composition 
	\[
		\nu=\nu\sub{p-1}\cdots\nu\sub1=id\sub{U}\times\L\sub{L\sub0}\times id\sub{\R{}}
	\]
	induces an unfolding
	\[
		U\times L\sub0\times\R{}\LARROW{\nu}U\times\widetilde{L\sub0}\times\R{}
	\]
	Now, by definition, $c=\widehat{c}\ \nu$. Applying \S\ref{def global thomather implies local thomather}-(b)
	to the stratified morphism $\widehat{\alpha}$ we get a stratified morphism $\widetilde{\alpha}$ between 0-length 
	stratified pseudomanifolds, from $U\times\widetilde{L\sub0}\times\R{}$ to $\X$, so 
	$\widetilde{\alpha}$ is smooth. We obtain a 
	commutative diagram
    \begin{center}
    \begin{picture}(70,70)(15,15)  \label{diagunbending}       
            \put(-50,10){$U\times c\left(L\sub0\right)$}
            \put(60,10){$X$}
            \put(-60,70){$U\times \widetilde{L\sub0}\times\mb{R}$}
            \put(60,70){$\X$}
 
            \put(-45,40){$^{c}$}
            \put(23,83){$_{\widetilde\alpha}$}
            \put(25,23){$_{\alpha}$}
            \put(70,40){$^{\L}$}
            \put(-35,65){{\vector(0,-1){40}}}
            \put(65,65){{\vector(0,-1){40}}}
            \put(5,75){\vector(1,0){40}}
            \put(5,15){\vector(1,0){40}}
        \end{picture}
    \end{center}\vskip2mm
	as desired.
\end{itemize}
\edem

\subsection{Smooth liftings}
We now study the smooth lifting of a Thom-Mather \linebreak morphisms.

\begin{Def}\label{def primary unfolding}
	The {\bf primary unfolding} of a given a finite length 
	Thom-Mather stratified pseudomanifold is the one we obtain with the
	iterated unbending process described on 
	\S\ref{Teo TM implies unfoldable}. 
\end{Def}

\begin{Pro}\label{pro thommather morphism lift into smooth map}
	Every Thom-Mather morphism between finite-length stratified Thom-Mather
	pseudomanifolds is unfoldable, in the sense of 
	\S\ref{def unfolding}-(4).
\end{Pro}
\bdem Let $X\ARROW{f}X'$ be a Thom-Mahter morphism.
According to \S\ref{Teo TM implies unfoldable} $X$ and $X'$
are unbendable and unfoldable.  By \S\ref{def global thomather implies local thomather}
$f$ is unbendable, $X\sub1=\widehat{X}$ and $X'\sub1=\widehat{X'}$ are Thom-Mather
stratified pseudomanifolds and the induced map 
\[
	X\sub1\ARROW{\hskip-2mm f\sub1=\widehat{f}}X'\sub1
\]
is a Thom-Mather morphism and satisfies $f\widehat{\L}=\widehat{\L'} f\sub1$.\vskip2mm 

After a finite number of iterated composition of unbendings, namely\linebreak  $n=\max\{l(X),l(X')\}$, 
we get two smooth primary unfoldings $\X=X\sub{n}\ARROW{\L}X$
and $\widetilde{X'}=X'\sub{n}\ARROW{\L'}X'$; we deduce that the respective iterated
$n$-th unbending $\widetilde{f}=f\sub{n}$ of $f$ satisfies $f\L=\L'\widetilde{f}$
and is smooth.
\edem

\begin{Teo}\label{Teo primary unfolding is a funtor}\hfill
	\begin{enumerate}
		\item The primary unfolding process is a functor in from the category of 
		finite length Thom-Mather stratified pseudomanifolds to the category
		of smooth manifolds.
		\item The primary unfolding of a Thom-Mather stratified pseudomanifold
		is unique up to Thom-Mather isomorphisms.
	\end{enumerate}
\end{Teo}
\bdem
	This is a consequence of \S\ref{pro thommather morphism lift into smooth map}.
\edem

\begin{center}
	{\bf Acknowledgments}
\end{center}
T. Guardia would like to thank the referee of {\sl Topology and its
Applications}, so as M. Paluzsny, for their useful and clever remarks. 
He also thanks the CDCH-UCV and FONACIT for some partial finnancial support of this research. 

\end{document}